\documentclass[letterpaper, 10 pt, conference]{ieeeconf}  

\IEEEoverridecommandlockouts                              
\usepackage{amsmath}
\usepackage{amsfonts}
\usepackage{amssymb}
\usepackage{color}
\usepackage{graphicx}
\usepackage{mathtools, balance}
\usepackage{xcolor} 
\usepackage{algorithm}
\usepackage[noend]{algpseudocode}
\usepackage{booktabs}
\usepackage{subcaption}
\usepackage[font=small,labelfont=bf]{caption}

\usepackage{tikz}
\usetikzlibrary{shapes,matrix,decorations.shapes}
\usetikzlibrary{arrows,decorations.pathmorphing,backgrounds,positioning,fit,matrix}
\usetikzlibrary{shapes,arrows,chains}
\usetikzlibrary[calc]
\usetikzlibrary[plotmarks]
\usetikzlibrary[patterns,decorations]

\newcommand{\R}{\mathbb{R}}

\newcommand{\psd}{\mathbb{S}}
\newcommand{\cs}{\mathcal{C}}

\newcommand{\gs}{\mathcal{G}}
\newcommand{\es}{\mathcal{E}}
\newcommand{\vs}{\mathcal{V}}

\DeclarePairedDelimiter{\abs}{\lvert}{\rvert}
\DeclarePairedDelimiter{\norm}{\lVert}{\rVert}

\DeclarePairedDelimiter{\Tr}{\text{Tr}(}{)}
\DeclarePairedDelimiter{\rank}{\text{rank}(}{)}

\DeclarePairedDelimiter{\vvec}{\text{vec}(}{)}

\DeclarePairedDelimiterX{\inp}[2]{\langle}{\rangle}{#1, #2}

\newtheorem{theorem}{Theorem}

\newtheorem{remark}{Remark}
\newtheorem{proposition}{Proposition}

\title{\LARGE \bf
Chordal Decomposition in Rank Minimized Semidefinite Programs with Applications to Subspace Clustering
}

\author{Jared Miller$^1$, Yang Zheng$^2$, Biel Roig-Solvas$^1$,  Mario Sznaier$^1$, Antonis Papachristodoulou$^3$
\thanks{$^1$J. Miller, B. Roig-Solvas and M. Sznaier are with the ECE Department, Northeastern University, Boston, MA 02115. (Emails: miller.jare@husky.neu.edu,  biel@ece.neu.edu, msznaier@coe.neu.edu).}
\thanks{$^2$Y. Zheng was with the Department of Engineering Science, University of Oxford. He is now with SEAS and CGBC at Harvard University, Cambridge, MA, 02138. (Email: zhengy@g.havard.edu).}
\thanks{$^3$A. Papachristodoulou is with the Department of Engineering Science, University of Oxford, Oxford, UK OX1 3PJ. (Email: antonis@eng.ox.ac.uk).}
\thanks{J. Miller, B. Roig-Solvas and M. Sznaier were partially supported by NSF grants  CNS--1646121, CMMI--1638234, IIS--1814631 and ECCS—1808381  and AFORS grant FA9550-19-1-0005.}
}

\begin{document}

\maketitle
\thispagestyle{empty}
\pagestyle{empty}

\begin{abstract}

Semidefinite programs (SDPs) often arise in relaxations of some NP-hard problems, and if the solution of the SDP obeys certain rank constraints, the relaxation will be tight. Decomposition methods based on chordal sparsity have already been applied to speed up the solution of sparse SDPs, but methods for dealing with rank constraints are underdeveloped. This paper leverages a minimum rank completion result to decompose the rank constraint on a single large matrix into multiple rank constraints on a set of smaller matrices. The re-weighted heuristic is used as a proxy for rank, and the specific form of the heuristic preserves the sparsity pattern between iterations. Implementations of rank-minimized SDPs through interior-point and first-order algorithms are discussed. The problem of subspace clustering is used to demonstrate the computational improvement of the proposed method.

\end{abstract}

\section{Introduction}

Semidefinite programs (SDPs) are a class of convex optimization problems that minimize a linear functional of a positive semidefinite (PSD) matrix under linear constraints. 
SDPs often arise as relaxations of some NP-hard problems, such as binary optimization, maxcut, and optimal power flow~\cite{vandenberghe1999applications}.
In particular, subspace clustering is an NP-hard problem that groups points originating from a union of subspaces and admits an SDP relaxation~\cite{cheng2016subspace}.  

In each case, SDP relaxations return equivalent optima if the solution satisfies rank constraints. Maxcut, optimal power flow, and subspace clustering all require rank-1 solutions.
%
The desire for a low-rank solution can be formulated as a rank-constrained SDP
\begin{equation} \label{Eq:PrimalSDP}
    \begin{aligned}
    \min_{X} \quad & \langle C,X \rangle \\
    \text{subject to} \quad & \langle A_i,X \rangle = b_i, i = 1, \ldots, m, \\
    & X \in \mathbb{S}^n_+,  \quad \text{rank}(X) \leq t,
    \end{aligned} 
\end{equation}
where $\inp{M}{N}=\Tr{M^T N}$ is an inner product and $C, A_1, \ldots, A_m \in \mathbb{S}^n$, $b\in\mathbb{R}^m$, and $t \in \mathbb{N}$ are problem data. Throughout this work, $\mathbb{R}^m$ is the $m$-dimensional Euclidean space, $\mathbb{S}^n$ is the space of $n \times  n$ symmetric matrices and $\mathbb{S}^n_+$is the subspace of symmetric PSD matrices.

While the rank-constrained SDP in \eqref{Eq:PrimalSDP} is in general NP-hard, a great deal of interest has been put in developing tractable rank proxies. 
Matrix factorization (\emph{i.e.}, $X = Y Y^T$)  upper bounds the possible rank by the width of $Y$, and Burer-Monteiro results may ensure global optimality after curvilinear optimization on the nonconvex low-rank manifold \cite{journee2010low}. 
%
Several convex relaxations of the rank constraint have been developed, one of the most popular being the nuclear norm $\norm{X}_* = \sum_{i} {\sigma_{i}(X)}$, and $\norm{X}_* = \Tr{X}$ if $X \in \mathbb{S}^n_+$ \cite{recht2010guaranteed}. The nuclear norm is the biconjugate of the rank function,  
and under  \emph{restricted isometry property} (RIP)/coherence conditions, the nuclear norm-minimized solution of an SDP is equivalent to the rank-minimal optimum \cite{candes2011tight}. RIP/coherence holds only in a very narrow set of problems  
since it weighs all singular values equally. 
Low rank solutions can be encouraged by penalizing smaller singular values.
The \emph{reweighted heuristic} is a linearization of log-det, and adds a penalty of $\Tr{W_t X}$ instead of the normal nuclear norm $\Tr{X}$, where $W_t = (X + \delta I)^{-1}$ updates at each iteration \cite{mohan2010reweighted}. 


The complexity of solving SDPs scales in a polynomial time w.r.t. the number of variables and constraints, and imposing rank penalties may further increase this complexity.
In sparse cases, only a small subset of entries of $X$ are used in the cost $C$ and constraints $A_i$, and the other entries guarantee $X \succeq 0$. 
When this sparsity pattern is (or can be extended to) a chordal graph, chordal decomposition theory can break up a large PSD cone ($X \succeq 0$) into a set of smaller coupled PSD cone ($X_k \succeq 0$)~\cite{vandenberghe2015chordal}.
Exploiting this structure reduces computational time, as the complexity of optimization problems is related to the graph's tree-width~\cite{chandrasekaran2012complexity}. Chordal decomposition of SDPs can effectively reduce dimension of sparse problems; see \emph{e.g.},~\cite{fukuda2001exploiting, zheng2017chordal} for details.

After performing the chordal decomposition and optimizing over the $X_k$'s, there exist multiple matrix completions to generate a valid $X \succeq 0$ from $X_k$. One such choice is the minimum rank completion \cite{dancis1992positive}, in which the minimal possible rank of the completion $X$ is the maximal rank among the blocks $X_k$. Numerical rounding on the eigenvalues of $X_k$ has already been used to reduce  of $\rank{X}$, but penalizing $\rank{X_k}$ in optimization was not considered~\cite{jiang2017minimum}. Minimum rank completions over linear matrix inequalities with general graphs has been performed in the context of optimal power flow, but few details were mentioned about how to penalize the rank of tree components~\cite{madani2017finding}. 

In this paper, we combine minimal rank completion and the reweighted heuristic to effectively solve large-scale rank-constrained SDPs. This heuristic penalty preserves the original problem's chordal sparsity while penalizing the rank of smaller matrices, and offers a rank guarantee based on the size of the maximum clique. 
We apply the proposed method to solve subspace clustering problems, which demonstrates the computational improvements.


The rest of this paper is structured as follows: Section \ref{Sec:preliminaries} introduces chordal graphs, sparse matrices, and the minimum rank completion problem. In Section \ref{Sec:Chordal}, we present an equivalent reformulation of~\eqref{Eq:PrimalSDP}. 
Implementations of chordal rank minimization in interior-point and first-order methods are discussed in Section \ref{Sec:algorithm}. These algorithms are used in Section \ref{Sec:subspace} to solve subspace clustering problems. Section~\ref{Sec:conclusions} concludes the paper. 
\section{Preliminaries} \label{Sec:preliminaries}

In this section, we cover chordal graphs, sparse (PSD) matrices, and the minimum rank completion. For a comprehensive treatment, the interested reader is referred to~\cite{vandenberghe2015chordal}.

\subsection{Chordal graphs}

An undirected graph $\gs(\vs, \es)$ is defined by a set of vertices $\vs=\{1,2,\ldots,n\} $ and edges $\es \subseteq \vs  \times \vs$. 
A cycle of length $N$ is a set of unique nodes ${v_k}$ such that $(v_1, v_2), \ldots, (v_i, v_{i+1}), \ldots,  (v_N, v_1) \in \es$. A chord is an edge that connects two nonconsecutive nodes in a cycle. An undirected graph is \emph{chordal} if all cycles of length no less than four have at least one chord~\cite{vandenberghe2015chordal}. 
A chordal extension  (a.k.a. completion) $\gs_c(\vs_c, \es_c)$ of graph $\gs(\vs, \es)$ is a chordal graph $\gs_c$ where $\vs = \vs_c$ and $\es \subseteq \es_c$. Finding a chordal extension with a minimal number of added edges is NP-hard, but efficient heuristics exist to give good chordal extensions~\cite{yannakakis1981computing}.

\begin{figure}[t]
	\centering
	\footnotesize
	\begin{tikzpicture}
	\matrix (m) [matrix of nodes,
	row sep = 0.8em,	
	column sep = 1em,	
	nodes={circle, draw=black}] at (-1.8,0)
	{ &  & & \\ 1 & 2 & 3 & 4 \\ & & &\\};
	\draw (m-2-1) -- (m-2-2);
	\draw (m-2-2) -- (m-2-3);
	\draw (m-2-3) -- (m-2-4);
	\node at (-1.8,-1.3) {(a)};
	
	\matrix (m3) [matrix of nodes,
	row sep = 0.8em,	
	column sep = 1.2em,	
	nodes={circle, draw=black}] at (1.8,0)
	{ & 2 & \\ 1 &  & 4 \\& 3 &\\};
	\draw (m3-1-2) -- (m3-2-1);
	\draw (m3-1-2) -- (m3-2-3);
	\draw (m3-2-1) -- (m3-3-2);
	\draw (m3-2-3) -- (m3-3-2);
	\draw (m3-1-2) -- (m3-3-2)[dashed];
	\node at (1.8,-1.3) {(b)};
	\end{tikzpicture}
	\caption{Examples of chordal graphs: (a) a path graph; (b) a triangulated graph (with dashed edge). Without the dashed edge there is a cycle of length 4 without a chord, so the graph is not chordal. The graph with the dashed edge is a chordal extension.}
	\label{F:ChordalGraph}
	\vspace{-2mm}
\end{figure}
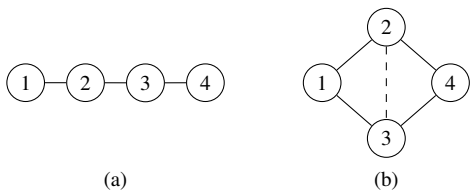

A clique $\cs$ with cardinality $\abs{\cs}$ is a subset of vertices in $\vs$ that form a complete subgraph. A maximal clique is a clique not contained inside another clique. Finding all max. cliques is NP-hard for general graphs, but can be computed on chordal graphs in linear time. Two  chordal graphs are shown in Fig.~\ref{F:ChordalGraph}, where the graph in Fig.~\ref{F:ChordalGraph}(a) has max. cliques $\mathcal{C}_i = \{i, i+1\}$, $i = 1, 2, 3$, and the graph in Fig.~\ref{F:ChordalGraph}(b) has max. cliques  $\mathcal{C}_1=\{1, 2, 3\}$ and $\mathcal{C}_2 = \{2, 3, 4\}$.

\subsection{Sparse matrices and chordal decomposition}
Considering the SDP~\eqref{Eq:PrimalSDP}, the graph $\gs(\vs, \es)$ arises from taking a union of the sparsity graphs corresponding to the data matrices $C, A_1, A_2, \ldots, A_m$. Precisely,
let $\hat{\es} = \es \cup \{(i, i), \forall i \in \vs\}$ be the edge set with self loops, and we define a cone of sparse symmetric matrices as:
    \[\psd^n(\es, 0) = \{X \in \psd^n \mid \ X_{ij} = 0, \ \forall (i, j) \not\in \hat{\es}\} \]
%
$\psd^n_+(\es, 0) = \psd(\es, 0) \cap \psd_+^n$ forms a cone of sparse PSD matrices. The dual cone $\psd^n_+(\es, 0)^* = \psd^n_+(\es, ?)$, which is the set of matrices that can be completed to be PSD with entries defined in $\es$ ($\exists M \succeq 0 \mid X_{ij} = M_{ij} \; \forall (i, j) \in \es^*) $.
Such a completion is not usually unique, as there may be multiple  $M$ associated to each $X$. For chordal graphs, $\psd^n_+(\es, ?)$ can be equivalently decomposed into a set
of smaller but coupled convex cones:

\begin{theorem} [Grone's Theorem {\cite{grone1984positive}}] \label{T:ChordalCompletionTheorem}
     Let $\mathcal{G}(\mathcal{V},\mathcal{E})$ be a chordal graph with a set of maximal cliques $\{\mathcal{C}_1,\mathcal{C}_2, \ldots, \mathcal{C}_p\}$. Then, $X\in\mathbb{S}^n_+(\mathcal{E},?)$ if and only if
    \end{theorem}
$$ E_{\mathcal{C}_k} X E_{\mathcal{C}_k}^T \in \mathbb{S}^{\vert \mathcal{C}_k \vert}_+,
    \qquad k=1,\,\ldots,\,p.$$
\vspace{-6mm}

In Theorem~\ref{T:ChordalCompletionTheorem}, $E_{\mathcal{C}_k}$ are $0/1$ entry selector matrices that index out components of $X$ involved in clique $\cs_k$. Grone's theorem provides a set equivalence $\psd^n_+(\es, ?) = \prod_{k=1}^{p}{\psd_{+}^{\abs{\cs^k}}}$ modulo overlaps between cliques, breaking a large PSD cone into a host of smaller PSD cones and equality constraints. This fact underpins the idea of much work that exploits sparsity in large-scale SDPs~\cite{vandenberghe2015chordal, fukuda2001exploiting, zheng2017chordal}.

\subsection{Minimum rank completion}

Given $X \in \psd_+(\es, ?)$, many choices of PSD completions are available after determining $X_k = E_{\mathcal{C}_k} X E_{\mathcal{C}_k}^T$, two of which are the maximum determinant completion and minimum rank completion. There exists a unique completion with maximum determinant with an explicit formula~\cite{barrett1989determinantal}. Minimum rank completions are not necessarily unique:

\begin{theorem}[Minimum rank completion~\cite{dancis1992positive}] \label{T:rankcompletion}
 Given a chordal graph $\mathcal{G}(\mathcal{V},\mathcal{E})$ with a set of maximal cliques $\{\mathcal{C}_1,\mathcal{C}_2, \ldots, \mathcal{C}_p\}$, for any $X\in\mathbb{S}^n_+(\mathcal{E},?)$, there exists at least one minimum rank PSD completion, where $ \text{rank}(X) = \max_{k} \text{rank}(E_{\mathcal{C}_k} X E_{\mathcal{C}_k}^T)$.
\end{theorem}

A procedure to perform minimum rank completion is Algorithm 3.1 in~\cite{jiang2017minimum}, which updates a factorization of the completion while proceeding through the elimination tree. 

We conclude this section with the following example:
\begin{figure}[!h]
    \centering
    \includegraphics[width=\linewidth]{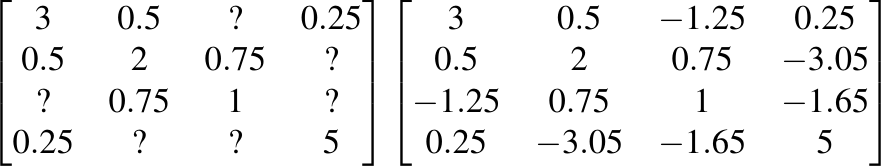}
    \caption{PSD Completable matrix (left) and rank-2 completion (right)}
    \label{fig:my_label}
\end{figure}

\noindent The maximum cliques of $\gs(\vs, \es)$ are $\{(1,2), (2,3), (1,4)\}$, and their induced submatrices are all PSD with rank 2, so there exists a completion with rank at least 2 (right matrix).

\section{Chordal decomposition in sparse SDPs with a rank constraint} \label{Sec:Chordal}

In this section, we first introduce a chordal decomposition approach in sparse SDPs with a rank constraint by combining Theorems~\ref{T:ChordalCompletionTheorem} and~\ref{T:rankcompletion}. Then, we discuss the application of reweighted heuristic  to the decomposed problem.

\subsection{An equivalent reformulation}
Let the rank constrained SDP~\eqref{Eq:PrimalSDP} be sparse with an \emph{aggregate sparsity pattern}  $\mathcal{G}(\mathcal{V}, \mathcal{E})$: \ 
$
    C \in \mathbb{S}^n(\mathcal{E},0), A_i \in \mathbb{S}^n(\mathcal{E},0), i = 1, \ldots, m.
$
Throughout this paper, we assume that $\gs(\vs, \es)$ is chordal (or has a suitable chordal extension) with a set of maximal cliques $\{\cs_k\}_{k=1}^p$.
Combining Grone's theorem (Theorem~\ref{T:ChordalCompletionTheorem}) with the minimum rank completion theorem (Theorem~\ref{T:rankcompletion}), 
leads to the following observation: 
\begin{proposition} \label{prop:decomposition}
    Suppose that Problem \eqref{Eq:PrimalSDP} is feasible and the problem data has an aggregate sparsity pattern $\mathcal{G}(\mathcal{V},\mathcal{E})$. Then,~\eqref{Eq:PrimalSDP} is equivalent to the following reformulation
%
    \begin{equation} \label{Eq:PrimalSDPdecomposition}
    \begin{aligned}
    \min_{X} \quad & \langle C,X \rangle \\
    \text{subject to} \quad & \langle A_i,X \rangle = b_i, i = 1, \ldots, m,\\
    & E_{\mathcal{C}_k} X E_{\mathcal{C}_k}^T\in \mathbb{S}^{\vert \mathcal{C}_k \vert}_+, \;\;\, k = 1, \ldots, p, \\
    & \text{rank}(E_{\mathcal{C}_k} X E_{\mathcal{C}_k}^T) \leq t, \, k = 1, \ldots, p,
    \end{aligned}
\end{equation}
    in the sense that~\eqref{Eq:PrimalSDP} and~\eqref{Eq:PrimalSDPdecomposition} have the same cost value, and their optimal solutions can be mutually recovered.
\end{proposition}

\begin{proof}
Thanks to the aggregate sparsity pattern, the cost function and the equality constraints in~\eqref{Eq:PrimalSDP} depend only on the elements $X_{ij}$ with $(i,j) \in \mathcal{E}^*$. The PSD constraint $X \in \mathbb{S}^n_+$ in~\eqref{Eq:PrimalSDP} can be equivalently replaced by a PSD completable constraint $X \in \mathbb{S}^n_+(\mathcal{E},?)$. 
The rest of proof directly follows the application of Theorems~\ref{T:ChordalCompletionTheorem} and~\ref{T:rankcompletion} to~\eqref{Eq:PrimalSDP}. We denote the optimal cost values to~\eqref{Eq:PrimalSDP} and~\eqref{Eq:PrimalSDPdecomposition} as $f^*_1$ and $f^*_2$ respectively.
    \begin{itemize}
      \item     First, assume $X^*_1$ is an optimal solution    to~\eqref{Eq:PrimalSDP} with cost $f^*_1 = \langle C,X^*_1 \rangle$ and $\text{rank}(X^*_1) \leq t$. Then, $X^*_1$ is also a feasible solution to~\eqref{Eq:PrimalSDPdecomposition}. Thus $f^*_2 \leq f^*_1$.

      \item Second, assume $X^*_2$ is an an optimal solution  to~\eqref{Eq:PrimalSDPdecomposition} with an cost value $f^*_2 = \langle C,X^*_2 \rangle$. According to Theorems~\ref{T:ChordalCompletionTheorem} and~\ref{T:rankcompletion}, we can find a PSD completion $\hat{X}^*_2$, where $\text{rank}(\hat{X}^*_2) = \max_{k} \text{rank}(E_{\mathcal{C}_k} X^*_2 E_{\mathcal{C}_k}^T) \leq t$. The PSD completion $\hat{X}^*_2$ is a feasible solution to ~\eqref{Eq:PrimalSDP}, indicating that $f^*_1 \leq f^*_2$.

    \end{itemize}

    Combining these facts, we know $f^*_1 = f^*_2$, and the optimal solutions to~\eqref{Eq:PrimalSDP} and~\eqref{Eq:PrimalSDPdecomposition} can be recovered from each other.
\end{proof}

One key feature of problem \eqref{Eq:PrimalSDPdecomposition} is that both the PSD and rank constraints are only imposed on multiple small symmetric matrices of smaller dimension rather than on the single large symmetric matrix. The minimum rank completion automatically yields an upper bound on the minimized full matrix rank according to maximum clique size.

\subsection{Rank relaxations}
\begin{figure}[t]
    \centering
    \includegraphics[width=0.95\linewidth]{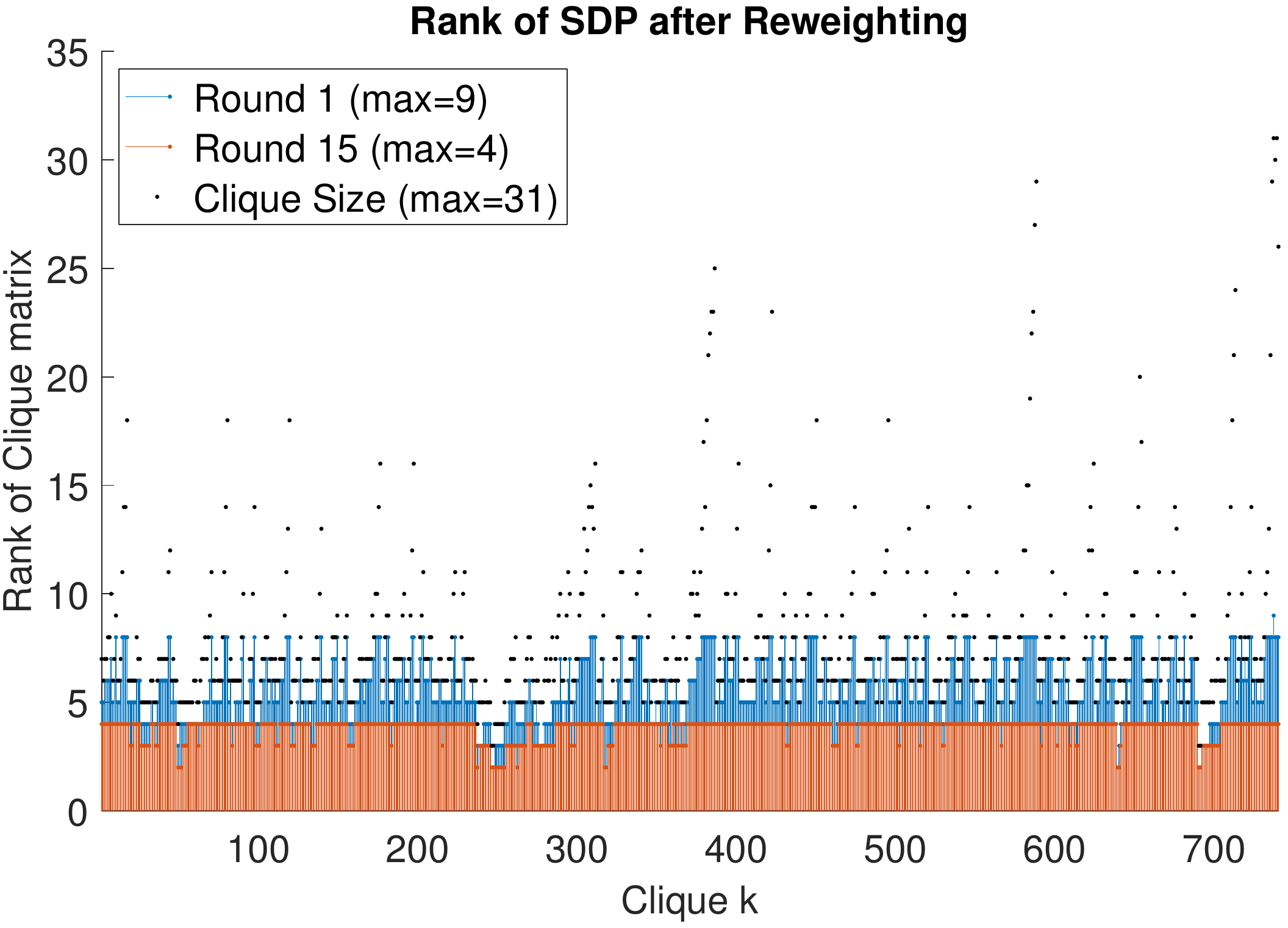}
    \caption{Example of reweighting heuristic (Algorithm 1) on a randomly generated Maxcut problem with 1000 vertices. There are 740 cliques with $\abs{\cs^{\text{max}}}=31$ (black dots). The maximum clique rank starts at $9$ (blue), and drops to $4$ (orange) after $15$ rounds of optimization.}
    \label{fig:maxcut_reweighting}
\end{figure}

In general, problem~\eqref{Eq:PrimalSDPdecomposition} is hard to solve due to the rank constraints. In this paper, we replace the hard rank constraint with a soft reweighted heuristic \cite{mohan2010reweighted}, leading to a standard SDP with chordal sparsity. The reweighted heursitic relaxes the hard $\rank{X} \leq t$ constraint to a soft $\inp{W}{X}$ term on the objective. For problem~\eqref{Eq:PrimalSDPdecomposition}, each rank constraint $\rank{E_{\cs_k} X E_{\cs_k}^T} \leq t$ is replaced by a soft penalty $\langle W_k, E_{\cs_k} X E_{\cs_k}^T\rangle $ for a reweighting matrix $W_k$ calculated in Algorithm 1. In general, the inner product between $W_k$ and $X_k$ promotes the concentration of energy of $X_k$ onto the dominant eigenspace of $X_{k,\text{old}}$ and thus incentivizes the reduction of its rank. Weights $W_k$ are normalized in Algorithm 1, and $\tau_k$ is a per clique regularization parameter. 
Reweighted heuristic is a local linearization of log-det and has been successfully used to solve rank minimization problems~\cite{mohan2010reweighted}.
Fig.~\ref{fig:maxcut_reweighting} demonstrates the rank reduction behavior of Algorithm 1 on a Maxcut SDP. 

The individual $W_k$'s can be combined together into an agglomerated penalty $W_\cs = \sum_{k=1}^pE_{\mathcal{C}_k}^TW_kE_{\mathcal{C}_k}$ by the cyclical property of Trace:
$$
    \begin{aligned}
          \sum_{k=1}^p\langle W_k,E_{\mathcal{C}_k} X E_{\mathcal{C}_k}^T \rangle 
        =&\left\langle \sum_{k=1}^pE_{\mathcal{C}_k}^TW_kE_{\mathcal{C}_k},X \right\rangle  \\
        =& \left\langle W_\cs,X \right\rangle. 
    \end{aligned} 
$$
This allows us to write the rank relaxed problem into the following form:
\begin{equation} \label{Eq:PrimalSDPdecompositionNuclears2}
    \begin{aligned}
    \min_{X} \quad & \left\langle C,X \right\rangle + \left\langle W_\cs,X \right\rangle \\
    \text{subject to} \quad & \langle A_i,X \rangle = b_i, i = 1, \ldots, m,\\
    & E_{\cs_k} X  E_{\cs_k}^T \succeq 0 , \forall\, k = 1, \ldots, p.
    \end{aligned}
    \end{equation}
The augmented cost matrix $C+ W_\cs \in \mathbb{S}^n(\mathcal{E},0)$ and $A_i \in \mathbb{S}^n(\mathcal{E},0)$, indicating that Problem~\eqref{Eq:PrimalSDPdecompositionNuclears2} has the same aggregate sparsity pattern as~\eqref{Eq:PrimalSDP}. If the original sparsity structure $\gs$ is non-chordal with a chordal extension $\gs_c(\vs, \es_c)$, the first reweighting iteration will have $W_\cs \in \psd^n(\es_c, 0)$. After the initial reweighting, there will be no more fill-in.


\begin{algorithm}[t] \label{alg:reweighted}
\caption{Chordal SDP with Reweighted Heuristic}
\begin{algorithmic}[1]
\Procedure{Chordal\_Rank}{}
\State $W_k \gets I \qquad \forall \, k = 1, \ldots, p$
\While {$X$ not converged}
\State $(X) \gets \text{optimum of \eqref{Eq:PrimalSDPdecompositionNuclears2} given}  W_k$
\State  $\tilde{W}_k \gets (E_{\cs_k} X_k E_{\cs_k}^T + \delta I)^{-1}$
\State  $W_k \gets \tau_k \tilde{W_k}/\norm{\tilde{W}_k}_2$
\EndWhile
\EndProcedure
\end{algorithmic}
\end{algorithm}

\begin{remark} \label{RK:reweighted}
As shown in Proposition~\ref{prop:decomposition}, Problems~\eqref{Eq:PrimalSDP} and~\eqref{Eq:PrimalSDPdecomposition} are equivalent. This tightness is lost when conducting rank relaxations. 
Applying the reweighted heuristic directly to Problem \eqref{Eq:PrimalSDP} requires the inverse of the big $X$ to update the weight $W = (X + \delta I)^{-1}$. As inverses are generically dense, the next iteration's cost $C + W$ would be dense and the sparsity pattern would be destroyed. 
In the proposed method, only the inverses of  the small $X_k$ matrices are required and the weight $W_\cs$ in \eqref{Eq:PrimalSDPdecompositionNuclears2} preserves the sparsity structure. Reweighted heuristic on problems \eqref{Eq:PrimalSDP} and \eqref{Eq:PrimalSDPdecomposition} produce different weights $W$ and $W_\cs$, and the resulting optima will not generally match as $W \neq W_\cs$. We note that using $W$ in the formulation of \eqref{Eq:PrimalSDP} will still penalize rank, but may settle at a high rank solution as the upper bound on clique size is lost. The convergence time of reweighted heuristic increases as the size of $X$ goes up, so in general it will take fewer iterations for a $W_\cs$-based scheme to converge than an algorithm using $W$ \cite{mohan2012iterative}.

\end{remark}
\section{Algorithm Implementations} \label{Sec:algorithm}



Problem~\eqref{Eq:PrimalSDPdecompositionNuclears2} is convex and readily solved by conic solvers, \emph{e.g.}, SeDuMi or Mosek. However, naively passing to solvers does not scale well to large-scale instances. In this section, we modify Problem~\eqref{Eq:PrimalSDPdecompositionNuclears2} to exploit its structure for adaptation in both interior point methods and first order methods. 

\subsection{Interior-point methods}

Interior point methods (IPM) such as SeDuMi will suffer from additional equality constraints introduced from the chordal decomposition forcing equality between cliques (clique-tree method). A different chordal extension $\tilde{\gs} \supset \gs$ may be used for optimization than the reweighting $\gs$ that forms $W_{\cs}$. SparseCoLo \cite{fujisawa2009user} or other methods may perform this decomposition to find an equivalent problem that is kinder to IPM, but reweighting is strictly based on the original pattern. IPM will be given Problem~\eqref{Eq:PrimalSDPdecompositionNuclears2} as input, after a possible decomposition.
%
%
%
%
In Algorithm 1, only the cost matrix $C + W_\cs$ changes between reweighting iterations, so no repeat conversions are necessary. 
Interior point methods generically find full rank solutions, so the low rank solution can be extracted by indexing out and rounding the cliques $X_k^*$ and then forming the minimal rank completion $X_r^*$.

\subsection{Alternating Direction Method of Multipliers (ADMM)}

Following~\cite{zheng2017chordal}, we build an ADMM algorithm to solve ~\eqref{Eq:PrimalSDPdecomposition}. Grone's theorem has a natural variable split $X_k = E_{\mathcal{C}_k} X E_{\mathcal{C}_k}^T$ that separates affine and PSD constraints. The resulting optimization problem is:
\begin{equation} \label{Eq:PrimalSDPdecompositionSplit}
    \begin{aligned}
    \min_{X,X_k} \quad & \langle C + W_{\cs},X \rangle \\
    \text{subject to} \quad & \langle A_i,X \rangle = b_i, i = 1, \ldots, m\\
    & X_k = E_{\mathcal{C}_k} X E_{\mathcal{C}_k}^T, \\
    & X_k \succeq 0 , \forall\, k = 1, \ldots, p.
    \end{aligned}
\end{equation}
Each iteration of ADMM has three steps~\cite{boyd2011distributed}: The affine constraints in $X$ are handled in step 1, the PSD projection in $X_k$ is form step 2, and the dual ascent on dual variables $\Lambda_k$ to enforce the variable-split is on step 3. First order algorithms such as ADMM sometimes converge slowly, but have a relatively low per-iteration cost.
Applying ADMM to~\eqref{Eq:PrimalSDPdecompositionSplit} requires solving subproblems at each iteration $t$ (with optional adjustment of $\rho$ between iterations following \cite{boyd2011distributed}):
\textbf{Step 1:} Solve the following quadratic program (QP)
  \begin{equation} 
    \begin{aligned}
    \min_{X} \: & \langle C+W_{\cs},X \rangle  +  \frac{\rho}{2}\sum_{k=1}^p 
    \left\|X_k^{(t)} - E_{\mathcal{C}_k} X E_{\mathcal{C}_k}^T + \frac{1}{\rho} \Lambda_k^{(t)}\right\|^2_F \\
    & \inp{A_i}{X} = b_i, \; i = 1, \ldots, m.
    \end{aligned}
\end{equation}
\textbf{Step 2:} Project onto $\psd_+^{\abs{\cs_k}}$ in parallel
  \begin{equation}
    \begin{aligned}
    \min_{X_k \succeq 0} \: &   \frac{\rho}{2}\left\|X_k - E_{\mathcal{C}_k} X^{(t+1)} E_{\mathcal{C}_k}^T + \frac{1}{\rho} \Lambda_k^{(t)}\right\|^2_F
    \end{aligned}
\end{equation}
\textbf{Step 3:} Update the multipliers by dual ascent
  \[  \Lambda_k^{(t+1)} = \Lambda_k^{(t)} + \rho\left(X_k^{(t+1)} - E_{\mathcal{C}_k} X^{(t+1)} E_{\mathcal{C}_k}^T\right). \]
The QP in step 1 can be solved by vectorizing all the variables. 
Entry selection on the clique $\cs_k$ can be replaced by a matrix $H_k$ where $H_k \vvec{X} = \vvec{E_{\mathcal{C}_k} X E_{\mathcal{C}_k}^T}$. Likewise $a_i = \vvec{A_i}$ is collated into $A$ and $b_i$ into $b$, $c = \vvec{C + W_{\cs}}$, and $v_k = \vvec{X_k^{(t)}  + \frac{1}{\rho} \Lambda_k^{(t)} }$. Step 1 is vectorized into:
  \begin{equation} 
    \begin{aligned}
    \min_{x} \quad & \inp{c}{x}  +  \frac{\rho}{2}\sum_{k=1}^p 
    \norm{H_k x - v_k}_2^2 \\
    \text{subject to} \quad &  A x = b.
    \end{aligned}
\end{equation}
The KKT system for the QP involves $x$ and dual variable $\omega$:
\[
\begin{bmatrix}
D & A^T \\
A & 0
\end{bmatrix}
\begin{bmatrix}
x \\ \omega
\end{bmatrix} = 
\begin{bmatrix}
\sum_{k=1}^p {H_k^T v_k}  - c \\ b
\end{bmatrix}.
\]
Since each $H_k$ is orthonormal, $H_k^T H_k$ is diagonal and $D = \sum_{k=1}^p {H_k^T H_k}$. The diagonal-offset lends itself nicely to block elimination and precomputed factorization; see~\cite{zheng2017chordal} for detailed discussions on a numerical implementation. 

Step 2 involves parallel PSD projections: find the eigendecomposition of $V_k = E_{\mathcal{C}_k} X^{(t+1)} E_{\mathcal{C}_k}^T - \frac{1}{\rho} \Lambda_k^{(t)}$ and keep the positive eigencomponents. 
If desired, the reweighting $W_{\cs}$ on $X$ can be moved to an individual weighting $W_k$ on $X_k$. 

\subsection{Homogeneous Self-Dual Embedding}
Problem~\eqref{Eq:PrimalSDPdecompositionSplit} can be solved through the homogeneous self-dual embedding (HSDE) framework~\cite{ye1994nl}. HSDE combines the primal and dual problems together to allow for the identification of infeasible SDPs. Solutions of chordally-sparse SDPs through HSDE have been already implemented~\cite{zheng2017hsde}, which is based on the general conic formulation introduced in \cite{o2016conic}. Each iteration of HSDE is comprised of a large block-sparse linear system, projection onto a product of multiple cones, and a dual update step (see~\cite{zheng2017hsde,o2016conic} for details). 
If $s$ is the concatenation of $\vvec{X_k}$ 
and $c = \vvec{C + W_{\cs}}$,
the linear system in HSDE $v=Qu$ is:
\[
\begin{bmatrix}
h \\ z \\ r \\ w \\ k
\end{bmatrix}
=
\begin{bmatrix}
& &  -A^T & -H^T & c\\
 & & & I & 0 \\
 A & & & & -b\\
 H & -I & & & \\
 -c^T & 0^T & b & &
\end{bmatrix}
\begin{bmatrix}
x \\ s \\ y \\ v \\ t
\end{bmatrix}.
\]

\section{Subspace Clustering} \label{Sec:subspace}

\begin{figure}[h]
    
    
\newlength{\twosubht}
    \newsavebox{\twosubbox}
\newlength{\firstfig}
    \newsavebox{\firstfigbox}
\newlength{\secondfig}
    \newsavebox{\secondfigbox}
    
    \sbox\twosubbox{%
      \resizebox{\dimexpr.95\linewidth}{!}{%
        \includegraphics[height=3.5cm]{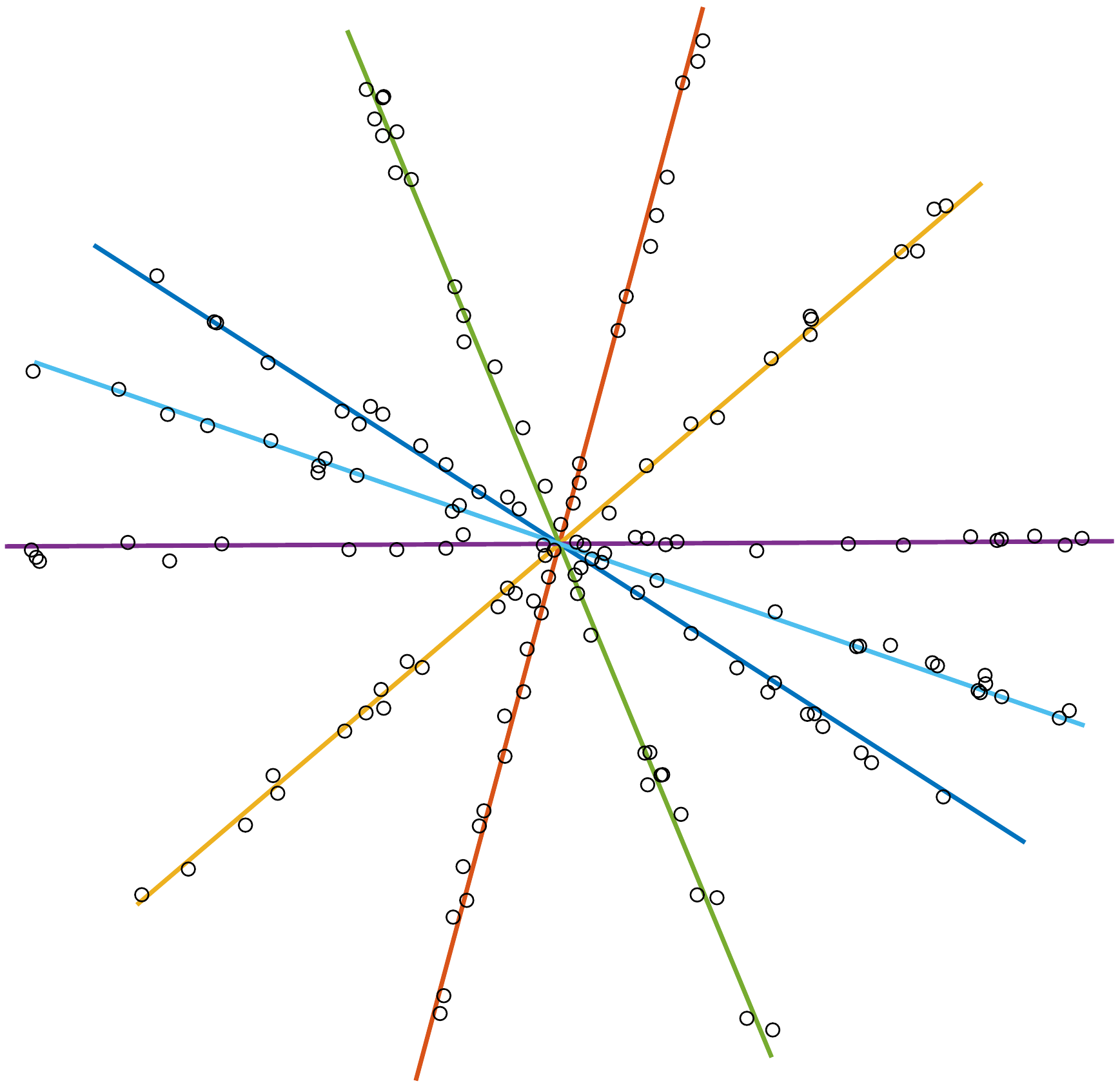}%
        \includegraphics[height=3.5cm]{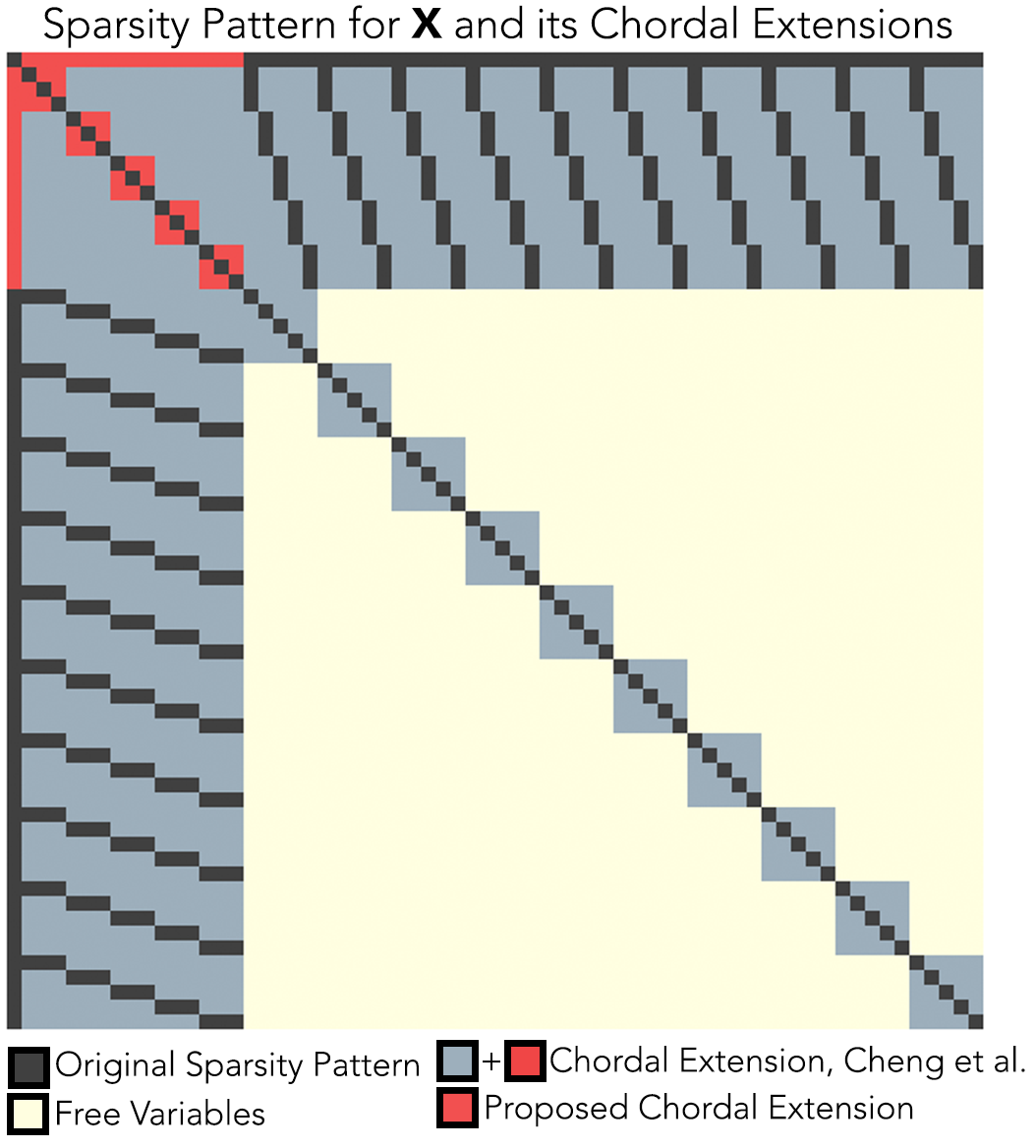}%
      }%
    }
    \setlength{\twosubht}{\ht\twosubbox}
    
    
    \centering

    \subcaptionbox{An application of SSC}{%
      \includegraphics[height=\twosubht]{img/SSC_successful.pdf}%
    }\quad
    \subcaptionbox{SSC sparsity patterns}{%
      \includegraphics[height=\twosubht]{img/SSC-sparsity.png}%

    }
    \caption{SSC example and sparsity patterns}
    \label{fig:SSCapplication}

\end{figure}

In this section, we demonstrate the performance of Algorithm 1 to solve subspace clustering (see Fig. \ref{fig:SSCapplication} a.). Subspace Clustering refers to the task of, given a set of $N_p$ points $x_{j}\in \mathbb{R}^D\;\forall_{j=1}^{N_p}$ sampled from $N_s$ subspaces, find the normals $r_i\in \mathbb{R}^D\;\forall_{i=1}^{N_s}$ of the sampled subspaces \cite{vidal2011subspace}. Alternatively, it can be posed as finding the set of binary labels $s_{i,j}\;\forall_{i=1}^{N_s}\forall_{j=1}^{N_p}$ that assign each point $x_j$ to the subspace spanned by the normal vector $r_i$. 
This problem arises in many practical applications including switched system identification from noisy input/output data (see \emph{e.g.}, \cite{paoletti2007identification,ChengWangSz14}), where each subspace normal is defined by the coefficients of each switching system.

In the general formulation, a point $x_j$ belongs to subspace $i$ if $r_i^T x_j = 0$. Under bounded noise, this orthogonality constraint is relaxed to $|r_i^T x_j|\leq \epsilon$, where $\epsilon$ is the noise bound. The task of subspace clustering can be cast as the following feasibility problem:
\begin{subequations}\label{SSC:main}
    \begin{align}
        \underset{r,s}{\textrm{find}} \qquad  & s_{i,j}|r_i^T x_j|\leq  s_{i,j}\epsilon,\quad \forall_{i=1}^{N_s},\forall_{j=1}^{N_p},\label{SSC:ortho}\\
        &s_{i,j} = s_{i,j}^2, \quad \forall_{i=1}^{N_s},\forall_{j=1}^{N_p},\label{SSC:binary} \\
        &\sum_{i=1}^{N_s} s_{i,j} = 1,\quad \forall_{j=1}^{N_p},\label{SSC:sum1} \\
        & r_i^Tr_i = 1, \quad \forall_{i=1}^{N_s},\label{SSC:norm1}
    \end{align}
\end{subequations}
where (\ref{SSC:ortho}) controls the orthogonality constraint and is only active when $s_{i,j}\neq 0$, \eqref{SSC:binary} enforces binary labels, \eqref{SSC:sum1} assigns every point $x_j$ to a subspace and \eqref{SSC:norm1} forces normal vectors to have unit norm (otherwise $r_i=0$ is a trivial feasible solution). 
Problem~\eqref{SSC:main} is nonconvex due to quadratic equality constraints and bilinear interactions between $s$ and $r$.

Cheng \emph{et. al. } proposed to reformulate \eqref{SSC:main} as an SDP by defining a matrix $X = [1,v]\,[1,v]^T$, with  
$$v = [r_1,\dots,r_{N_s},s_{1,1},\dots,s_{N_s,1},\dots,s_{N_s,N_j}]$$ containing all the variables of \eqref{SSC:main} in vectorized form \cite{cheng2016subspace}. $X$ is a symmetric PSD matrix of size $(1 + N_s(D+N_p))$, and all the constraints in \eqref{SSC:main} become linear with respect to the entries of $X$ at the cost of adding a non-convex rank 1 constraint. The reweighted nuclear norm heuristic was then employed to relax the rank constraint into a convex problem.
\footnote{Due to the particular structure of the subspace clustering problem, enforcing a rank 1 constraint on a particular principal submatrix of $X$ is equivalent to enforcing rank 1 on the overall matrix. The interested reader is referred to Appendix A of~\cite{cheng2016subspace} for the proof of this exact relaxation.}



We note that problem \eqref{SSC:main} lacks bilinear interactions between different $s_{i,j}$ terms or any interactions between $s_{\bar{i},j}$ and $r_{\hat{i}}$ when $\bar{i}\neq \hat{i}$. As a result, only a very reduced number of entries of $X$ are actually used in the constraints of \eqref{SSC:main}, leading to a very sparse pattern. The sparsity pattern for an SSC problem with parameters $(D=3, N_p=10, N_s=5)$ is shown in dark grey in Figure \ref{fig:SSCapplication}.b. To exploit this underlying sparsity, \cite{cheng2016subspace} proposed to solve \eqref{SSC:main} with variables in a chordal extension, shown by the union of light grey and red cells in Figure \ref{fig:SSCapplication}.b.
In this paper, we propose a reduced chordal extension of $X$, shown in red in Figure \ref{fig:SSCapplication}.b. The new extension is linear in $N_s$ and $N_p$, and imposes $N_s$ $[D+1]$-size rank-1 constraints instead of a single larger $[1+N_s D]$ constraint. A summary of the clique matrix sizes $\abs{\cs_k}$ is:

\begin{figure}[!h]
    \centering
    \includegraphics[width=\linewidth]{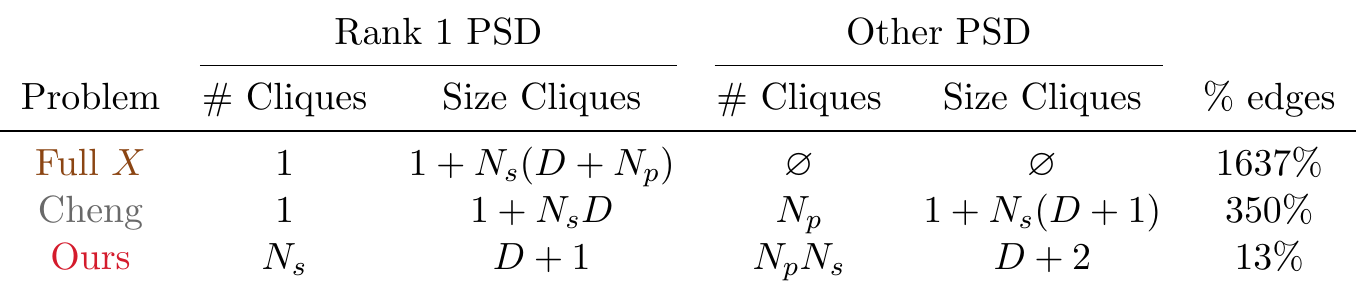}
    \captionof{table}{Clique sizes in SSC. \%edges measures size of chordal extension over {\color[RGB]{60,60,60}baseline $\es$} (variables in (\ref{SSC:main})).  Parameters are $(D=3, N_p=10, N_s=5$) as in Fig. \ref{fig:SSCapplication} b. }
    \label{fig:SSC_extension_sizes}
\end{figure}

\begin{figure}[!h]
    \centering
    \includegraphics[width=\linewidth]{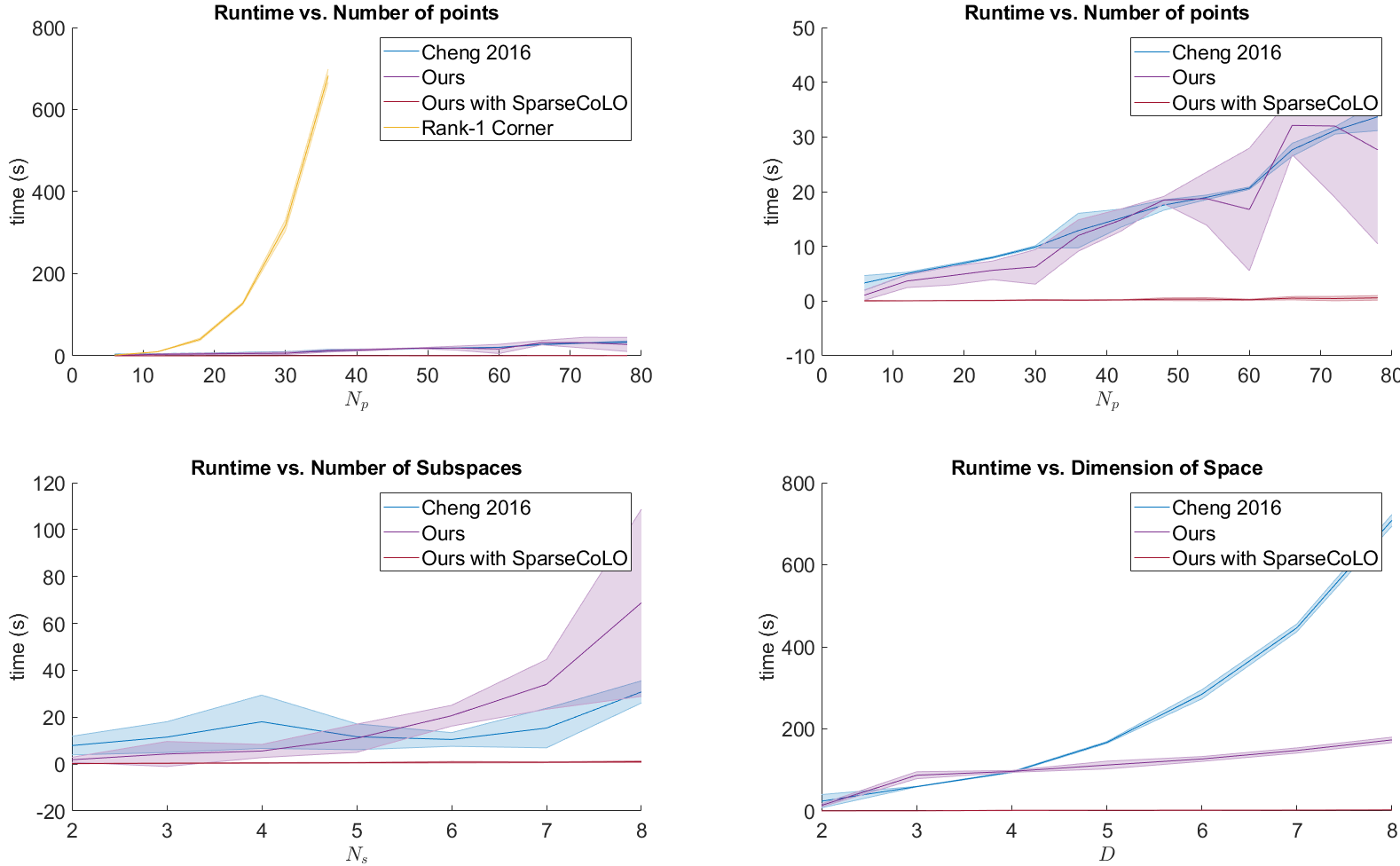}
    \caption{SSC Experiments, changing individually $N_p$ (top), $N_s$ (bottom left), and $D$ (bottom right).}
    \label{fig:SSCexperiments}
\end{figure}

For our extension, each clique $v_{i,j} = [1, r_i, s_{i,j}]$  where $r_i$ includes all coordinates of the normal. Experiments were run on SSC using YALMIP and Mosek. Data was generated by randomly choosing $N_s$ subspace normals $r_i \in \R^D$, then sampling $N_p$ total points $x_j$ from these subspaces. Points were corrupted by uniform noise with bounds $\pm\epsilon$. $W_{\cs}$-based regularization is applied until a rank-1 solution ($\sigma_1/(\sum{\sigma_i}) > 99\%$) is found or until 20 iterations. Results are recorded in Fig. \ref{fig:SSCexperiments}, all plots are mean$\pm$1 stdev. 

\begin{figure}[!h]
    \centering
    \includegraphics[width=\linewidth]{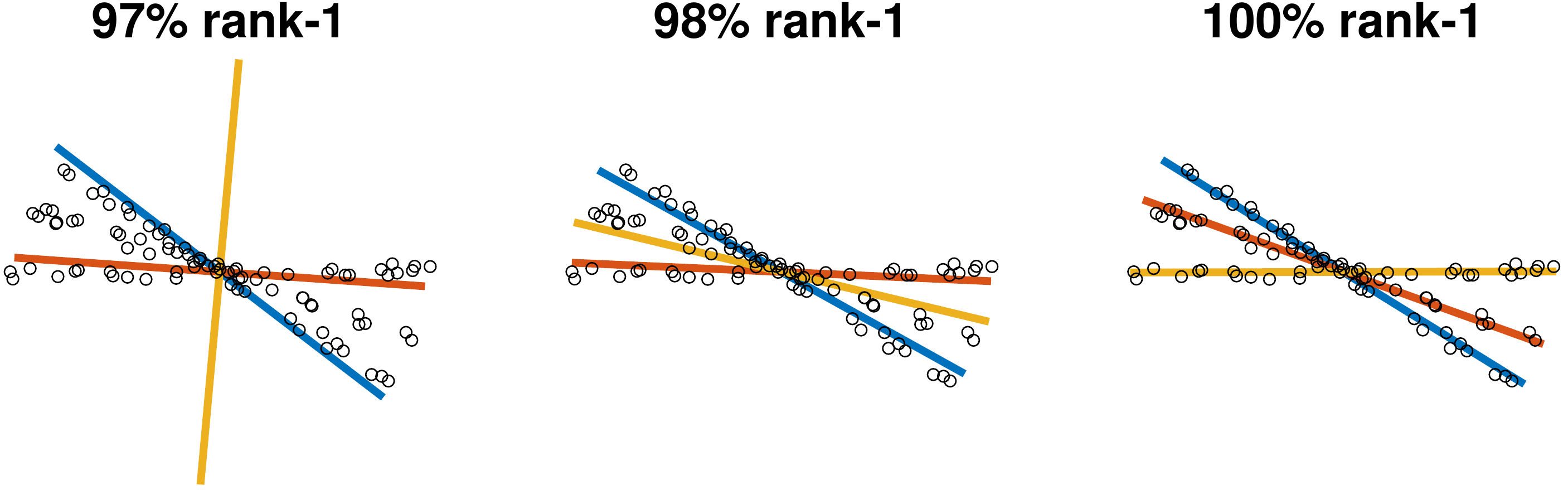}
    \caption{SSC SDP convergence to satisfying Problem~\eqref{SSC:main} as solution approaches rank-1. $N_s = 3, N_p = 90, D = 2, 
    \epsilon = 0.15$}
    \label{fig:SSCRank}
\end{figure}

Experiments are run on five penalty schemes: imposing that $X$ is rank-1, that the top-left corner of $X$ is rank-1, with the Cheng scheme, by  Problem~\eqref{Eq:PrimalSDPdecompositionNuclears2}, and by preprocessing Problem~\eqref{Eq:PrimalSDPdecompositionNuclears2} with SparseCoLO. The SparseCoLO scheme uses $W_{\cs}$ for rank penalization, and finds another decomposition before passing to Mosek. As an example with $N_s=4, N_p = 80, D=5$, Problem~\eqref{Eq:PrimalSDPdecompositionNuclears2} has 320 PSD blocks of size $D+2=7$, with 320 rank-1 constraints of size $D+1=6$. For this problem, SparseCoLO generates 25 PSD blocks with size between 12 and 30, and $W_{\cs}$ regularization and rank-guarantees are performed on $6 \times 6$ sub-blocks. In all tests, Problem~\eqref{Eq:PrimalSDPdecompositionNuclears2} with SparseCoLO is an order of magnitude faster than other problem formulations.

The top left plot of Fig. 
\ref{fig:SSCexperiments} has $N_s=3, N_p = 6*(1\ldots13), D = 3$. Imposing $\rank{X}=1$ directly never converges to a rank-1 solution. Setting $\rank{X_{\text{top left corner}}}=1$ is shown in yellow, and with $N_p=42$ each run takes 2 hours (not shown). The top right plot zooms in, and shows the Cheng formulation, Problem~\eqref{Eq:PrimalSDPdecompositionNuclears2}, and Problem~\eqref{Eq:PrimalSDPdecompositionNuclears2} as preprocessed by SparseCoLO \cite{fujisawa2009user}. Adding in clique-overlap constraints in Problem~\eqref{Eq:PrimalSDPdecompositionNuclears2} degrades performance as compared to Cheng, but with SparseCoLO preprocessing it wins out. In the bottom left,  $N_s=2\ldots8, N_p = 120, D = 2$ over 12 tests. $\rank{X}=1$ and $\rank{X_\text{corner}}=1$ do not converge to rank-1 solutions, and in this case Problem~\eqref{Eq:PrimalSDPdecompositionNuclears2} without SparseCoLO is slower than Cheng. The bottom right plot has $N_s=4, N_p = 80, D = 2\ldots8$ over 8 tests, and as $D$ increases our proposed formulation clearly beats Cheng. Variability in the plots is often due to iterations needed to find a low rank solutions. In some cases, rank-1 solutions are found in the first iterations, others require 20 iterations to become nearly rank-1 ($95\%$). These experiments also involve all $\tau_k=1$. Unequal weights of cliques (possibly by size of block) is of future interest.

\section{CONCLUSIONS} \label{Sec:conclusions}




In this paper, we combined the minimum rank completion with reweighted heuristic to solve rank-minimization problems with chordal sparsity. We showed that rank constraints only need to be placed on the cliques $X_k$, and that using $W_{\cs}$ maintains the sparsity pattern. We discussed implementations of chordal rank minimized SDPs by interior-point and first-order methods. We demonstrated the scalability and efficiency of the chordal decomposition for rank minimization in the specific example of subspace clustering. We expect that these gains will hold in many other chordally sparse rank-minimized SDPs, such as group (e.g. $\mathbb{Z}_2$) synchronization. Future work includes utilizing other rank surrogate functions and applying chordal rank minimization to more general polynomial and rational optimization problems.

\balance


\bibliographystyle{IEEEtran}
\bibliography{chordal_rank_paper}

\end{document}